\documentclass[12pt]{amsart}
\usepackage{amssymb}
\usepackage[all]{xy}

 \rm

\renewcommand{\)}{\right)}
\renewcommand{\[}{\left[}
\renewcommand{\]}{\right]}
\newcommand{\<}{\langle}
\renewcommand{\>}{\rangle}
\newcommand{\x}{\times}
\renewcommand{\bar}{\overline}

\newcommand{\st}{\:|\:}

\newcommand{\N}{{\mathbb{N}}}
\newcommand{\C}{{\mathbb{C}}}

\renewcommand{\phi}{\varphi}

\renewcommand{\H}{{\mathcal{H}}}
\newcommand{\BH}{{\mathcal{B}}(\H)}
\newcommand{\A}{{\mathcal{A}}}

\theoremstyle{plain}
\newtheorem{thm}{Theorem}

\newtheorem{cor}[thm]{Corollary}

\theoremstyle{definition}

\theoremstyle{remark}

\title{Inductive algebras for  compact groups}

\author{Promod~Sharma}
\author{M.~K.~Vemuri}
\address{Department of Mathematical Sciences\\ IIT (BHU)\\
         Varanasi 221 005\\ INDIA}

\begin{document}

\begin{abstract}
Inductive algebras for a compact group are self-adjoint.
\end{abstract}


\maketitle
\thispagestyle{empty}

\section{Introduction}\label{S:intro}
Let $G$ be a separable locally compact group and $\pi$ an irreducible
unitary representation of $G$ on a separable Hilbert space $\H$.  Let
$\BH$ denote the algebra of bounded operators on $\H$.  An {\em
inductive algebra} is a weakly closed abelian subalgebra $\A$ of $\BH$
that is normalized by $\pi(G)$, i.e., $\pi(g) \A \pi(g)^{-1} = \A$ for each
$g \in G$.  If we wish to emphasize the
dependence on $\pi$, we will use the term $\pi$-inductive algebra.
A {\em maximal inductive algebra} is a maximal element of the set of
inductive algebras, partially ordered by inclusion.

The identification of inductive algebras can shed light on the
possible realizations of $\H$ as a space of sections of a homogeneous
vector bundle (see e.g.\ \cite{sharma2022, Stegel-2004,Stegel-2006,sl2r,rcr}).
For self-adjoint maximal inductive algebras there is a precise result
known as Mackey's Imprimitivity Theorem, as explained in the
introduction to \cite{Stegel-2004}.  Inductive algebras have also found
applications in operator theory (see e.g.\ \cite{Koranyi-2014, iahs}).

In \cite{Raghavan-2005}, it was shown that finite dimensional
inductive algebras for a connected group are trivial.  However,
the title of \cite{Raghavan-2005} is somewhat misleading, as finite
groups can have non-trivial finite dimensional inductive algebras.

In this note we show that inductive algebras for a compact group are
self-adjoint.  This is significant because, in general, the
classification of self-adjoint inductive algebras is easier than the
classification of all inductive algebras.  This is because the methods
of spectral theory are available only in the former case.  Also, unlike
in the classification work cited above, we do not need to assume maximality.

In Section \ref{S:L_infty} we prove some results about subalgebras of
$L^{\infty}(X,\mu)$, which will be used in the proof of our main theorem,
but which are also of independent interest.

\section{Subalgebras of $L^\infty$}\label{S:L_infty}

\begin{thm}\label{T:sa}
Let $(X, \mu)$ be a measure space.  The algebra $L^{\infty}(X,\mu)$ is
finite dimensional if and only if all of its subalgebras are self-adjoint.
\end{thm}

\begin{proof}
Assume first that $L^{\infty}(X,\mu)$ is finite dimensional.  Observe that
under this hypothesis, if $f\in L^{\infty}(X,\mu)$, then there exists
a simple function $s$ such that $f=s$ almost everywhere (see
\cite[Prop 3.4.2]{Cohn} and \cite[\S13.3, Cor.\ 6]{Royden}).

Let $\A \subseteq L^{\infty}(X,\mu)$ be a subalgebra. Let
$\{f_1,f_2,\dots,f_n\}$ be a basis for $\A$, and choose simple
functions $s_1,s_2,\dots,s_n$ such that $f_j=s_j$ (a.e.),
$j=1,2,\dots,n$.  Define a  map $\mathbf{s}: X\to \C^n$ by
\begin{equation*}
\mathbf{s}(x)=\(s_1(x),s_2(x),\dots,s_n(x)\).
\end{equation*}
Since simple functions attain only finitely many values, $\mathbf{s}(X)$ is
finite, and we may write
\begin{equation*}
\mathbf{s}(X)\setminus\{0\}=\{\mathbf{v}_1,\mathbf{v}_2,\dots,\mathbf{v}_m\},
\end{equation*}
for some $m\in\N$.

Put $A_0=\mathbf{s}^{-1}(0)$ and
\begin{equation*}
A_k = \mathbf{s}^{-1}(\mathbf{v}_k), \qquad k=1,2,\dots,m.
\end{equation*}

Then $\{A_k\}_{k=0}^m$ are disjoint, and
\begin{equation*}
X=\bigcup_{k=0}^m A_k.
\end{equation*}
Let $v_{kj}$ denote the $j$-th component of the vector $\mathbf{v}_k$.
Observe that if $x\in A_k$ then $s_j(x)=v_{kj}$, $j=1, \dots, n$,
$k=1, \dots, m$, i.e., each $s_j$ is constant on each $A_k$.
Therefore
\begin{equation*}
s_j \in \operatorname{span}\{\chi_{A_k}\}_{k=1}^m, \quad j=1, \dots, n.
\end{equation*}
Therefore $\A \subseteq \operatorname{span}\{\chi_{A_k}\}_{k=1}^{m}$.

Fix distinct $h,k \in \{1, \dots, m\}$.  Since
$\mathbf{v}_h \ne \mathbf{v}_k$, there exists $j=j(h,k)$ such that
$v_{hj}\ne v_{kj}$.
Since $\mathbf{v}_h\ne 0$, there exists $l=l(h)$ such that $v_{hl}\ne 0$.
Observe that
\begin{equation*}
\phi_{hk}=(s_j - v_{kj}) s_l \in \A.
\end{equation*}

If $x \in A_k$, then
\begin{equation*}
\begin{aligned}
\phi_{hk}(x)
=&\; (s_j(x) - v_{kj}) s_l(x)\\
=&\; (v_{kj} - v_{kj}) v_{kl}\\
=&\; 0,
\end{aligned}
\end{equation*}
and if $x\in A_h$, then
\begin{equation*}
\begin{aligned}
\phi_{hk}(x)
=&\; (s_j(x) - v_{kj}) s_l(x)\\
=&\; (v_{hj} - v_{kj}) s_l(x)\\
=&\; (v_{hj}-v_{kj}) v_{hl}\\
\neq\;& 0.
\end{aligned}
\end{equation*}
Put
\begin{equation*}
\psi_{hk}=\frac{\phi_{hk}}{(v_{hj}-v_{kj}) v_{hl}}.
\end{equation*}
Then $\psi_{hk} \in \A$, $\psi(x)=1$ if $x\in A_h$
and $\psi(x)=0$ if $x\in A_k$.

Since $\chi_{A_h} = \prod_{k\neq h} \psi_{hk}$, it follows that
$\chi_{A_h} \in \A$, $h=1, \dots, m$.
Therefore $\A = \operatorname{span}\{\chi_{A_k}\}_{k=1}^{m}$.
Therefore $\A$ is self-adjoint.

Assume now that $L^\infty(X,\mu)$ is infinite dimensional.  We claim first
that $X$ has a sequence $\{E_n\}_{n=1}^\infty$ of disjoint measurable subsets of
positive measure.  Indeed, there exists a real valued function
$f\in L^{\infty}(X,\mu)$ such that
\begin{equation*}
-\infty < \operatorname{ess\, inf} f < \operatorname{ess\, sup} f < \infty.
\end{equation*}
Put
\begin{equation*}
c=\frac{\operatorname{ess\, inf} f + \operatorname{ess\, sup} f}{2}.
\end{equation*}
Let
\begin{equation*}
Y= \{x\in X|f(x)>c\}, \quad\text{and}\quad
Z= \{x\in X|f(x)\leq c\}.
\end{equation*}
Then $Y$ and $Z$ are disjoint measurable sets, and by the definitions of
essential supremum and essential infimum $\mu(Y)>0$ and $\mu(Z)>0$.  Since
\begin{equation*}
L^{\infty}(X)\cong L^{\infty}(Y)\oplus L^{\infty}(Z),
\end{equation*}
either $L^{\infty}(Y)$ or $L^{\infty}(Z)$ must be infinite
dimensional, say $\dim L^{\infty}(Z)=\infty$.  Let $E_1=Y$.
We may iterate the previous argument with $Z$ in place of $X$
to produce the required sequence.

Choose points $e_n \in E_n$, and let $\A$ consist of all $f\in L^\infty(X,\mu)$
which are constant on each $E_n$, $n=1,2, \dots$, and
\begin{equation*}
\lim_{m\to\infty} \frac{f(e_{2m+1})-f(e_1)}{(1/m)} =
 i \lim_{m\to\infty} \frac{f(e_{2m})-f(e_1)}{(1/m)}.
\end{equation*}
It is easy to check that $\A$ is a subalgebra of $L^\infty(X,\mu)$.
Now, define $f:X\to\C$ by
\begin{equation*}
f(x)=
\begin{cases}
0           & \text{if $x\in E_1$},\\
\frac{1}{n} & \text{if $x\in E_n$, $n$ even},\\
\frac{i}{n} & \text{if $x\in E_n$, $n>1$ odd}.
\end{cases}
\end{equation*}
Then $f\in\A$ but $\bar{f}\notin\A$.  Therefore $\A$ is not self-adjoint.
\end{proof}

\section{Compact Groups}

\begin{thm}\label{T:compact}
Let $G$ be a compact group and $\pi$ an irreducible unitary representation of
$G$ on a Hilbert space $\H$.  If $\A \subseteq \BH$ is a $\pi$-inductive
algebra, then $\A$ is self-adjoint.
\end{thm}

\begin{proof}
By the Peter-Weyl theorem, $\H$ is finite dimensional.

Let $\mathcal{N}$ denote the set of nilpotent elements in $\A$ (the
nilradical of $\A$).  Let
\begin{equation*}
\mathcal{K}=\{x \in \H \st Tx=0, \quad \forall T \in \mathcal{N}\}.
\end{equation*}
By (a trivial case of) Engel's theorem \cite[\S3.3]{Humphreys-1972},
$\mathcal{K}\ne 0$.  Observe that $\mathcal{N}$ is normalized by
$\pi(G)$, so $\mathcal{K}$ is $\pi(G)$-invariant.  However, since
$\pi$ is irreducible, it follows that $\mathcal{K}=\H$, whence
$\mathcal{N}=0$.

Let $\A^*$ denote the space of linear functionals on $\A$.
For each $\lambda \in \A^*$, let
\begin{equation*}
\H_\lambda =\{v\in \H \st Tv=\lambda(T)v \quad \text{for all $T\in\A$}\},
\end{equation*}
and
\begin{equation*}
\Lambda = \{\lambda\in \A^*\st \H_\lambda \neq 0 \}.
\end{equation*}
Then $\Lambda$ is a finite set.

Since $\A$ is abelian, and $\mathcal{N}=0$, the Jordan-Chevalley decomposition
\cite[\S4.2]{Humphreys-1972} implies that
\begin{equation}\label{E:JordanChevalley}
\H=\bigoplus_{\lambda\in \Lambda} \H_{\lambda}.
\end{equation}

Let $\<\cdot,\cdot\>$ denote the inner product of $\H$.
There exists an inner product $\<\cdot,\cdot\>_1$ on $\H$ such that
$\H_\lambda$ and $\H_\mu$ are orthogonal with respect to $\<\cdot,\cdot\>_1$
if $\lambda\ne\mu$.  Let $\sigma$ denote the Haar probability measure on
the compact group $G$.  By Schur's lemma (see \cite{brocker}), there exists a
constant $c$ such that
\begin{equation*}
\<v,w\> = c \int_G \<\pi(g)v,\pi(g)w\>_1\, d\sigma.
\end{equation*}

If $g\in G$ and $\lambda\in \A^*$, define $g\lambda: \A \to \C$ by
\begin{equation*}
g\lambda(T)=\lambda(\pi(g)^{-1}T\pi(g)), \quad T\in\A.
\end{equation*}
This defines an action of $G$ on $\A^*$, which preserves $\Lambda$.

Note that for any $g \in G$, $\lambda\in\A^*$ and $v\in \H_\lambda$,
$\pi(g)v\in \H_{g\lambda}$.
Also $\lambda \neq \mu$ implies $g\lambda \ne g\mu$.  Therefore, if
$\lambda \neq \mu$, $v\in \H_\lambda$ and $w \in \H_\mu,$ then
\begin{equation*}
\begin{aligned}
\<v,w\>
=\;& c \int_{G}\<\pi(g)v,\pi(g)w\>_1\, d\mu\\
=\;& 0.
\end{aligned}
\end{equation*}

Therefore $\H_\lambda$ and $\H_\mu$ are orthogonal with respect to
$\<\cdot,\cdot\>$ if $\lambda\ne\mu$.

Observe that if $\lambda\in\Lambda$, then $\lambda$ is multiplicative.
Indeed, if $\lambda\in\Lambda$, then there exists
$v\in\H_\lambda\setminus\{0\}$.  Therefore, if $T_1, T_2 \in \A$, then
\begin{equation*}
\lambda(T_1T_2)v=T_1T_2v=T_1(\lambda(T_2)v)=\lambda(T_2)T_1v=
\lambda(T_2)\lambda(T_1)v.
\end{equation*}
Since $v\ne 0$, it follows that $\lambda(T_1T_2)=\lambda(T_1)\lambda(T_2)$.

It follows that the map $\mathcal{G}:\A\to L^{\infty}(\Lambda)$
(with respect to counting measure) defined by
\begin{equation*}
\[\mathcal{G}(T)\](\lambda)=\lambda(T).
\end{equation*}
is an algebra homomorphism.

Since $\Lambda$ is finite, $L^\infty(\Lambda)$ is finite dimensional, and so
$\mathcal{G}(\A)$ is self-adjoint by Theorem \ref{T:sa}.

Let $T\in\A$.  Then there exists $T_1\in \A$ such that
$\mathcal{G}(T_1)=\bar{\mathcal{G}(T)}$,
i.e., $\lambda(T_1)=\bar{\lambda(T)}$ for all $\lambda\in\Lambda$.  We
claim that $T_1=T^*$, i.e., that
\begin{equation*}
\<Tv,w\>=\<v,T_1w\> \quad\text{for all $v,w\in\H$}.
\end{equation*}
By (\ref{E:JordanChevalley}) it suffices to check this assuming that
$v\in \H_\lambda$ and $w\in \H_\mu$ for $\lambda, \mu \in \Lambda$.

If $\lambda=\mu$, then
\begin{equation*}
\begin{aligned}
\<Tv, w\>
=&\; \<\lambda(T)v, w\>\\
=&\; \<v, \bar{\lambda(T)}w\>\\
=&\; \<v, T_1w\>.
\end{aligned}
\end{equation*}

If $\lambda\ne\mu$, then $\<v,w\>=0$, and so
\begin{equation*}
\begin{aligned}
\<Tv, w\>
=&\; \<\lambda(T)v, w\>\\
=&\; 0, \quad\text{and}\\
\<v, T_1w\>
=&\; \<v, \mu(T_1)w\>\\
=&\; 0.
\end{aligned}
\end{equation*}
\end{proof}

\begin{cor}\label{C:finite}
Let $G$ be a finite group and $\pi$ an irreducible unitary representation of
$G$ on a Hilbert space $\H$.  If $\A \subseteq \BH$ is a $\pi$-inductive
algebra, then $\A$ is self-adjoint.
\end{cor}

In view of Raghavan's theorem \cite{Raghavan-2005}, it might appear that
Corollary \ref{C:finite} may be used whenever Theorem \ref{T:compact}
is applicable.  However, that is not the case.  Indeed, if
$G=\mathrm{O}(2)$, the group of orthogonal $2\x 2$ matrices, then $G$
is compact and not abelian, but its group of components is abelian.  If
$\pi$ is an irreducible representation of $G$ of dimension greater than
one, then Theorem \ref{T:compact} implies that all $\pi$-inductive algebras
are self-adjoint, but Corollary \ref{C:finite} is not applicable.

\bibliographystyle{amsplain}
\bibliography{v9-iafcg.bib}
\end{document}